\documentstyle[amscd]{amsart}

\newcommand{\U}{{U}}

\newcommand{\bfZ}{{\bf  Z}}
\newcommand{\cB}{{\cal  B}}
\newcommand{\cZ}{{\cal  Z}}
\newcommand{\cH}{{\cal  H}}

\newsymbol\ltimes 226E
\newsymbol\rtimes 226F
\newsymbol\boxtimes 1202
\newsymbol\twoheadrightarrow 1310
\newsymbol\subsetneq 2328

\newcommand{\sur}{\twoheadrightarrow}

\renewcommand{\P}{{\cal P}}

\newcommand{\N}{{\cal N}}

\newcommand{\bu}{\bullet}

\newcommand{\A}{{\cal A}}
\newcommand{\B}{{\cal B}}

\newcommand{\F}{{\cal F}}
\newcommand{\Ftil}{\tilde{\cal F}}
\newcommand{\G}{{\cal G}}
\newcommand{\Gtil}{\tilde{\cal G}}
\newcommand{\ftil}{\tilde f}
\newcommand{\jtil}{\tilde j}
\newcommand{\tilF}{\tilde{\cal F}}

\newcommand{\Dpp}{{D^{p,\geq 0}}}

\newcommand{\Dpm}{{D^{p,\leq 0}}}

\newcommand{\cupl}{{\bigcup\limits}}

\newcommand{\xbar}{{\overline x}}
\newcommand{\pbar}{{\overline p}}

\newcommand{\Obar}{{\overline O}}
\newcommand{\Ubar}{{\overline U}}

\newcommand{\I}{{\cal I}}

\newcommand{\Zet}{{\Bbb Z}}
\newcommand{\D}{{\Bbb D}}
\newcommand{\DS}{{\Bbb DC}}
\newcommand{\DC}{{\Bbb DC}}

\newcommand{\Mtop}{X^{top}}
\newcommand{\Ztop}{Z^{top}}
\newcommand{\Utop}{U^{top}}

\newcommand{\Xtop}{X^{top}}

\newcommand{\uHom}{{\underline {Hom}}}

\renewcommand{\O}{{\cal O}}

\newcommand{\Xtil}{{\tilde X}}

\newcommand{\iso}{{\widetilde \longrightarrow}}

\newcommand{\imbed}{\hookrightarrow}

\newcommand{\jbar}{{\overline j}}

\newcommand{\proof}{{\it Proof\ \ }}
\newcommand{\epf}{\square}

\newcommand{\tauleqn}{\tau_{\leq n}^{stand}}
\newcommand{\taugeqn}{\tau_{\geq n}^{stand}}

\newtheorem{Thm}{Theorem}
\newtheorem{Cor}{Corollary}
\newtheorem{Lem}{Lemma}
\newtheorem{Prop}{Proposition}

\theoremstyle{definition}
\newtheorem{Def}{Definition}

\theoremstyle{remark}

\newtheorem{Rem}{Remark}
\newtheorem{Ex}{Example}

\def\square{\hbox{\vrule\vbox{\hrule\phantom{o}\hrule}\vrule}}

\begin{document}

\title[]{Perverse Coherent Sheaves (after Deligne)}
\author{Roman Bezrukavnikov}


\begin{abstract}
This note is mostly
an exposition of an unpublished result of Deligne [D], which
introduces an analogue of perverse $t$-structure [BBD] on the derived category
of coherent sheaves on a Noetherian scheme with a dualizing 
complex.  Construction extends to the category
of coherent sheaves equivariant under an action of an algebraic group;
though proof of the general statement in this case 
 does not require  new ideas, it provides examples (such as sheaves
on the nilpotent cone of a semi-simple group equivariant under the adjoint
 action) where construction of  coherent ``intersection cohomology''
sheaves works. 
\end{abstract}

\maketitle

\section{Introduction} 
Let $X$ be a reasonable stratified topological space; or let
$X$ be a reasonable scheme, stratified by locally closed subschemes.
Let $D$ be the full subcategory in, respectively, derived category
of sheaves on $X$, or in  the derived category of etale sheaves  
on $X$,  consisting of complexes smooth along the stratification.

For an integer-valued function $p$ (perversity) on the set of strata
Beilinson, Bernstein and Deligne \cite{BBD}
defined a  $t$-structure on the category $D$; the objects of corresponding
abelian category (core of the $t$-structure) are called perverse sheaves.

The question addressed in this note is whether an analogous construction
can be carried out for the derived category of coherent sheaves on a reasonable
scheme. Surprisingly, the answer is positive (with some modifications),
easy, and not widely known (though was known to Deligne for a long time, 
see \cite{D}). 

Let us summarize the difference
between  the coherent case considered here, and
the constructible
case treated in \cite{BBD}.

First, in the coherent case
we can not work with  complexes ``smooth'' along a 
given stratification, for the corresponding subcategory in $D^b(Coh)$ is not
a full triangulated subcategory. (If $f$ is a function whose divisor intersects
the open stratum, then the cone of the morphism $\O \overset{f}{\to} \O$
has singularity on the open stratum). This forces us to define perversity as a
function on the set of generic points of all 
irreducible subvarieties, i.e. on the
topological space of a scheme. 

The second, more essential
 difference is that  in the derived category of coherent sheaves
the functor $j^*$ of  pull-back  under an open imbedding $j$
does not have adjoint functors. 
Recall that in  constructible situation
the right adjoint to $j^*$ is the functor  $j_*$ of 
  direct image, and the left adjoint is the functor  $j_!$
of extension by zero.
In coherent set-up the functor $j_*$ is defined in the larger
 category of quasi-coherent sheaves (Ind-coherent sheaves),
 while $j_!$ is defined in the Grothendieck dual
 category introduced in Deligne's appendix to \cite{H} (consisting 
of Pro-coherent sheaves).

It turns out, however, that in the proof 
 of the existence of perverse $t$-structure one can use
instead of the object  $j_!(\F)$ (where $j:U\imbed X$ is an open imbedding)
any extension $\Ftil$ of $\F$ to $X$, such that
the restriction  of $\tilF$ to $X-U$ has no cohomology above certain degree
(depending on the perversity function).
 If the perversity function is
{\it monotone}  (see Definition \ref{mon} below) it is very easy to construct
such $\Ftil$. Applying the Grothendieck-Serre duality
to this construction, we get a substitute for $j_*(\F)$, which exists
if the perversity function is comonotone.
 Otherwise the proof is parallel to that in \cite{BBD}.

Thus the
 $t$-structure is constructed not for an
 arbitrary perversity function, but only for a monotone
and comonotone one. (In the topological situation one also needs this
condition to get a $t$-structure on the whole derived category
of constructible sheaves, rather than on the category corresponding to a fixed
stratification.)

In \cite{D} Deligne used the
Grothendieck's Finiteness Theorem (\cite{SGA2}, VIII.2.1) to show that
 the formulas for 
$\tau ^p_{\leq 0}$, $\tau^p_{\geq 0}$ of \cite{BBD}, a'priori making sense
in a larger category containing $D^b(Coh)$, give in fact objects of $D^b(Coh)$,
provided the perversity function is monotone and comonotone
 (see also Remark \ref{Grfin}).

\medskip

The results on the existence of a ``perverse'' $t$-structure carry
over to the case of $G$-equivariant coherent sheaves, where $G$ is a
 (reasonable)  algebraic group acting on a (reasonable) scheme.
In this case  perversity $p(x )$ must be assigned only to points
$x $ of the scheme, which are invariant under the connected component of
identity of $G$, as an equivariant sheaf is anyway ``smooth along
the orbits.''

Although the general formalism for the equivariant category is very
similar to the non-equivariant one (to the extent that we found it easier
not to treat the two cases separately),
 there is one construction which works
in the equivariant case only.
Namely, the definition of the minimal (Goresky-MacPherson, or IC)
 extension functor
$j_{!*}$ works only when the perversity function is {\it strictly}
monotone and comonotone. Though formally the proof of this statement
works both
in the equivariant and non-equivariant (=equivariant with trivial $G$)
situations, the statement can be nonempty in the equivariant case only.
Indeed, it is easy to see, that a strictly monotone and comonotone
perversity function
 exists only if $G$ acts on the scheme with finite number of orbits,
and dimensions of two adjacent orbits differ at least by two.
If this is the case, an obvious analogue
of the usual desciption of irreducible perverse
sheaves as minimal extensions of local systems is valid,
and the core of the $t$-structure is Artinian (in contrast with
the core of the standard $t$-structure). An example of this situation
is provided by the nilpotent cone of a semi-simple algebraic group
over a field of characteristic zero, equipped with the adjoint action
(see Remarks at the end of the note).

The exposition would probably look better (and work in greater generality)
if the notion of a stack was used; however, my ignorance confined me
to the language of equivariant sheaves (rather than the equivalent 
language of sheaves on the quotient stack).

\medskip
It should be
 quite clear from the above that this paper does not contain original
results of the author.

{\bf Acknowledgements.} 
I am much obliged to  Pierre Deligne for valuable explanations, and for a
kind permission to use his unpublished results.

I thank Dima Arinkin, Alexander Beilinson, Victor Ginzburg,
and Dmitry Panyushev for discussions and references.

Quite separately, I express my deepest gratitude to Leonid Positselskii:
without his participation in the early stages of the work it might never have 
been done.

This work was started in the spring of 1999 when the author was a
member at the Institute for Advanced Study. I thank IAS for excellent work
conditions, and NSF grant DMS 97-29992 for financial support.

\section{Preliminaries}\label{secprel}
In this section we collect some standard Lemmas needed in the exposition.
The reader familiar with basic algebraic geometry certainly will not need
our proofs.

Let $X$ be a scheme over a base scheme $S$; we
  denote the category of  coherent 
(respectively, quasi-coherent) sheaves on $X$ by $Coh_X$, $QuasiCoh_X$
(or simply $Coh$, $QuasiCoh$ if confusion is not likely).
Let  $G$ be an affine  group scheme over
$S$, acting on $X$.

We will assume that $S$ is Noetherian, $X$, $G$ are of finite type over $S$,
and $S$ admits a dualizing complex (in the sense of 
\cite{H} \S V.2);
the structure morphism $f_G:G\to S$ is assumed to be flat of finite type
 and Gorenstein
(i.e. $f_G^!(\O_S)$ is locally free). 


 The category of $G$-equivariant coherent 
(respectively, quasi-coherent) sheaves on $X$ is denoted by
 $Coh_X^G$, $QuasiCoh_X^G$.
The forgetful functor $Forg:QuasiCoh^G\to QuasiCoh$
 has the right adjoint $Av: 
\F \mapsto a_* pr^* \F$, where $pr:G\times X\to X$ and 
$a:G\times X\to X$ are respectively the projection and the action.
(Here $Av$ stands for ``averaging''.)
Since $G$ is affine and flat, $Av$ is exact, and
  the canonical morphism $\F\to Av(\F)$ is an 
injection. Hence $QuasiCoh^G$ has enough injectives, because
 $QuasiCoh$ does.

\begin{Lem} \label{union}
 Any $G$ equivariant quasi-coherent sheaf $\F$ on $X$ is the union of its
$G$-equivariant coherent subsheaves.
\end{Lem}

\proof Let $a:\F\to \O(G) \otimes_{\O(S)} \F$ denote the coaction.
For a coherent (possibly non-equivariant)
subsheaf $\F_0 \subset \F$ let $\F_0^G\subset \F$ 
be the preimage under $a$ of $\O(G) \otimes_{\O(S)} \F_0$.
Then one readily checks that $\F_0^G$ is an equivariant coherent subsheaf,
and that $\cupl_{\F_0} \F_0^G=\F$. \epf

\begin{Cor} \label{odinhren}            
 For $?=b$,  or $-$ the category
$D^?(Coh^G)$ is equivalent to the full 
subcategory of $D^?(QuasiCoh^G)$ consisting of complexes with coherent 
cohomology. 

\end{Cor}

  \proof It suffices to check that for a bounded above complex $\F^\bu$ of
equivariant quasicoherent sheaves, whose cohomology is coherent, 
 the set of quasiisomorphic  equivariant coherent subcomplexes in $\F^\bu  $
  is nonempty and filtered under inclusion; and that any equivariant coherent
  subcomplex in $\F^\bu$ lies in an equivariant coherent quasiisomorphic 
      subcomplex.
This follows from  Lemma \ref{union} by a standard argument. Namely,
let $\cZ^i, \cB ^i\subset \F^i$ denote, respectively, the kernel and the 
image of the differential. We construct by descending induction in $i$ a
 coherent equivariant subsheaf 
$\F^i_c \subset \F^i$ satisfying the two properties: $d(\F^i_c)=
\F^{i+1}_c\cap \cB^{i+1}$; and $(\F^i_c\cap \cZ^i)\sur \cZ^i/\cB^i=\cH^i$.
If we are given a coherent subcomplex $\G^\bu\subset \F^\bu$
we can choose $\F_c^i $ to satisfy also $\F_c^i\supset \G^i$. \epf

We will denote the full subcategories in $D^+(QuasiCoh)$,
$D^+(QuasiCoh^G)$ consisting of complexes with coherent cohomology by
$D^+_c(QuasiCoh)$, $D^+_c(QuasiCoh^G)$ respectively.

\begin{Cor}\label{exte} Let                  
$U\subset X$
be an open $G$-invariant subscheme. 
For any $\F,\G\in D^b(Coh^G(U))$ and a 
morphism $f:\F\to \G$ 
there exist $\Ftil, \Gtil\in D^b(Coh^G(X))$, and a morphism
$\ftil:\Ftil\to \Gtil$, such that $\ftil|_U \cong f$. For any two
$\Ftil', \Ftil''\in D^b(Coh^G(X))$  and an isomorphism
 $f:\Ftil'|_U
\cong \Ftil''|_U$ there exists $\Ftil\in D^b(Coh^G(X))$, and morphisms
$f':\Ftil'\to \Ftil$, $f'':\Ftil''\to \Ftil$ such that 
$f'|_U$, $f''|_U$ are isomorphisms, and $f''|_U\circ f =f'|_U$.
\end{Cor}

\proof Let $\F^\bu$ be a finite complex of equivariant quasicoherent sheaves
on $X$, such that cohomology of $C^\bu|_U$ is coherent. 
A construction similar to the one used in the proof of Corollary
\ref{odinhren} shows that the set of coherent equivariant subcomplexes
$\F_c^\bu \subset \F^\bu$, 
such that imbedding of $\F_c^\bu|_U \imbed \F^\bu|_U$ is quasiisomorphism,
is nonempty and filtered under inclusions. Moreover, any  equivariant coherent
  subcomplex in $\F^\bu$ lies in such a subcomplex.
The statement follows. \epf

We will write $\Xtop$ for the topological space of a scheme $X$.
If $x \in \Mtop$ is a point of $X$ (respectively $Z\subset X$ is a locally 
closed subscheme), then $i_x :\{pt\}\imbed \Xtop$ (respectively 
 $i_Z:Z \imbed X$) will denote the imbedding.

 We will use the same notation for a functor on
an abelian category and its derived functor. In particular, 
for $x\in \Xtop$  the  functors
$i_x^*:D^?(Coh)\to D^?(O_x-mod)$, and $i_x^!:D^+(QuasiCoh)\to D^+(O_x-mod)$
are derived of  respectively an exact, and of a left exact functor. 
The functor $i_X^!$ factors through the derived category of
torsion $\O_x$ modules, and has finite homological dimension
(because so does the functor $j_*$ where $j$ is an open imbedding in a
Noetherian scheme).

{\it 
Everywhere below we will assume that $Coh$ has enough locally free objects.
Also, dealing with equivariant categories, we will assume that $Coh^G$
has enough locally free objects.}

Then both in $D(Coh)$ and in $D(Coh^G)$
 internal Hom (denoted by $\uHom$) can
be computed as derived functor in either of the two variables, and
commutes with the forgetful functor from the equivariant to the nonequivariant
 category; also
for a ($G$-equivariant) morphism $f:Z\to X$ the 
coherent pull-back functor $f^*$ 
is defined in both categories and commutes with  forgetful functor.

\begin{Lem}\label{equiv} Let    
  $Z\subset X$ be a locally closed  ($G$-invariant)
 subscheme,   and
 $n$ be an integer. Let $x \in \Xtop$ be a generic point of $Z$. Then
  
a) For $\F \in D^-(Coh^G)$ we have $i_x ^*(\F )\in D^{\leq n}(\O _x -mod)$ iff
there  exists an open ($G$-invariant) subscheme $Z^0\subset Z$, $Z^0\owns x$,
 such that $i_{Z^0}^*(\F )\in D^{\leq n} (Coh_{Z_0}^G)$;

b)   For $\F \in D^+_c(QuasiCoh^G)$ 
we have $i_x ^!(\F )\in D^{\geq n}(\O_x -mod)$ 
iff there
  exists  an open  ($G$-invariant) subscheme $Z^0\subset Z$, $Z^0\owns x$,
such that $i_{Z^0}^!(\F )\in D^{\geq n} (Coh_{Z_0}^G)$.
\end{Lem}

\proof Existence of an open ($G$-invariant) subscheme $Z^0\subset Z$
as in (a) is equivalent to $i_x^*i_Z^*(\F)\in D^{\leq n}((\O_Z)_x-mod)$.
(Indeed, if the last equality holds, then we can let $Z^0$ be 
the complement in $Z$ to support of ${\cal H}^k (i_Z^*(\F))$, $k>n$;
the converse is obvious.)

We can rewrite $i_x^*i_Z^*(\F)
=i_x^*(\F)\overset{L}{\otimes}_{\O_x} \O(Z)_x$. Since the functor
of tensor product with $ \O(Z)_x$ over $\O_x$ is right exact, and kills no
finitely generated $\O_x$ modules by the Nakayama Lemma, we see that 
the top cohomology of $i_x^*(\F)\overset{L}{\otimes}_{\O_x} \O(Z)_x$
and of $i_x^*(\F)$ occur in the same degree. This proves (a).

Similarly,  the second condition in (b) says that $i_x^!
i_Z^!(\F)=i_x^*i_Z^!(\F)\in  D^{\geq n}((\O_Z)_x-mod)$
(the equality here is, of course, due to the fact $x$ is generic in $Z$).
We rewrite
$i_x^!i_Z^!(\F)
=RHom_{\O_x}( \O(Z)_x, i_x^!(\F))$, and see that 
 the lowest cohomology of  $i_x^!(\F)$ and of 
$RHom_{\O_x}( \O(Z)_x, i_x^!(\F))$ occur at the same degree, because
$Hom_{\O_x}( \O(Z)_x,\underline{\ \ })$ is left exact, and kills no 
torsion module, while cohomology of $i_x^!(\F)$ is a torsion $\O_x$-module.
\epf

\begin{Lem}\label{lim}
Let $i:\bfZ \imbed \Xtop$ be imbedding of a closed
$G$-invariant subspace. 

a) (cf. e.g.
\cite{H}, Theorem V.4.1)
For any $\F \in D^-(Coh^G)$, $\G\in D^+(QuasiCoh^G)$ we have
$$                                 
Hom(\F ,i_*i^!(\G))=\varinjlim_Z Hom(\F ,i_{Z*}i_Z^!(\G)),
$$                                 
where $Z$ runs over the set of closed $G$-invariant
subschemes of $X$ with the underlying
topological space $\bfZ$.

b) If $\F\in D^b(Coh^G)$ is such that the cohomology sheaves
${\cal H}^i(\F)$ are supported on $\bfZ$, then there exists
a closed $G$-invariant subscheme $Z\subset X$, $\Ztop=\bfZ$, such
that $\F\cong i_{Z*} (\F_Z)$ for some $\F_Z\in Coh^G(Z)$. 
\end{Lem}

\proof a)
Let us represent $\F $ by a bounded above complex $P_\F $ of locally 
free coherent equivariant sheaves, and $\G$ by a bounded below complex $I_\G$
 of injective 
quasicoherent equivariant sheaves.
If $\I$ is an injective object of $QuasiCoh^G(X)$, and $Z\subset X$
is a closed $G$-invariant subscheme, then $i_Z^!(\I)\in QuasiCoh^G(Z) $
is injective; hence locally free equivariant sheaves on $X$ are adjusted
to $Hom(\underline{\ \ },i_{Z*}i_Z^!(\I))$.
Thus  $RHom(\F, i_{Z*}i_Z^!(\G))=RHom(i_Z^*(\F), i_Z^!(\G))$ is computed
by the complex $Hom^\bullet(P_\F, i_{Z*}i_Z^!(I_\G))$. On the other hand, 
$RHom(\F ,i_*i^!(\G))$ is computed by the complex 
$Hom^\bullet(P_\F, i_* i^!(I_\G))$, as $ i_* i!$ sends injective
objects of $QuasiCoh^G$ into injective ones. We have
$ i_*i^!(I_\G)=\cupl_Z i_{Z*}i_Z^!(I_\G)$; and also 
$Hom^\bullet
(\P_\F, i_*i^!(I_\G))=\cupl_Z Hom^\bullet (P_\F ,i_{Z*}i_Z^!(I_\G))$,
 because $P_\F$ is a bounded above complex of coherent sheaves.
This implies the Lemma. 

b) The category $QuasiCoh^G_\bfZ(X)$ of $G$-equivariant quasi-coherent sheaves 
supported on $Z$ has enough injectives; moreover they are also injective
as objects of the larger category $QuasiCoh^G(X)$ (this follows from the
corresponding statement for non-equivariant sheaves, since $Av$
preserves sheaves supported on $\bfZ$). Hence
 (see e.g  \cite{H}, Proposition I.4.8) $\F$ is quasiisomorphic to a
finite   complex
of quasicoherent sheaves supported on $\bfZ$. As in the proof
of Corollary \ref{odinhren} we can represent this complex
 as a union of quasiisomorphic
equivariant  coherent subcomplexes; any such subcomplex 
is  supported on a closed
subscheme $Z$, $\Ztop\subset \bfZ$.
\epf

\begin{Def}  An equivariant dualizing complex on $X$ is an object
$\DS ^G\in D^b(Coh^G)$, such that every $\F\in D^b(Coh^G)$ is $\DS$-reflexive,
i.e. 
the natural transformation
$\F \to \uHom(\uHom(\F,\DS),\DS)$ is an isomorphism.
\end{Def}

\begin{Lem} $\F\in D^b(Coh^G)$ is an equivariant dualizing complex iff 
$Forg(\F)$ is a dualizing complex.
\end{Lem} 

\proof The 'if' direction is clear because $\uHom$ commutes with
the forgetful functor. The 'only if' follows from \cite{H}, Proposition V.2.1,
which says, in particular, that if the structure sheaf $\O$ is $\DS$ reflexive,
then $\DS$ is a dualizing complex. Since $\O$ obviously lies in the image of 
the forgetful functor, we see that $Forg(\DS^G)$ is a dualizing complex.
\epf

\begin{Prop}\label{dual_exist} In the above
assumptions
 $X$ admits an equivariant dualizing complex. 
\end{Prop}

\proof  According to \cite{BBD}, Theorem 3.2.4 an object of the derived
category of sheaves on a cite can be given locally provided negative
local Ext's from the object to itself vanish. Applying it to the
covering $G\times X\to X$ in
the  cite of flat $G$-schemes over $X$, we see that it is enough to provide
an isomorphism $\pi^*(\DC)=a^*(\DC)$ (here that $\pi:G\times X\to X$, and 
$a:G\times X\to X$ are the projection and the action maps), satisfying
an  associativity constraint on $G\times G\times X$. 
Since $f_G:G\to S$ is Gorenstein, the sheaf $f_G^!(\O_S)$ is invertible;
 the group structure on $G$ provides then a canonical
isomorphism $f_G^*=f_G^!$
(as follows e.g. from Remark in \cite{H}, pp 143-144).
 Hence $\pi^*(\DC)=a^*(\DC)$ are both
canonically isomorphic to $f_{G\times X}^!(\DC_S)$, which provides the desired 
isomorphism. The associativity constraint follows from functorial properties
of $f^!$.
\epf

\begin{Rem} Suppose that we make an additional assumption
that the structure morphism $X\to B$ is {\it
 equivariantly embeddable}, i.e. can be presented as a composition 
 $X\overset{\iota}{\imbed} \Xtil\to B$, where $\Xtil$ is
a smooth $B$-scheme with a $G$-action, and $\iota$ is a $G$-equivariant closed 
imbedding (the Sumihiro embedding Theorem
\cite{Su} (see also \cite{KK})
 guarantees that this assumption is satisfied if $S$ is the spectrum of
an algebraically closed field of characteristic 0, and $X$ is a normal
quasiprojective variety). Then 
the Proposition becomes evident, for we can set $\DS_X^G\overset{\rm def}{=}
\iota^!(\Omega^{top}_{\Xtil})$, the definition of $\iota^!$ for a closed 
imbedding being straightforward.
\end{Rem}

\section{Perverse coherent sheaves}\label{s1}
\subsection{Construction of the $t$-structure}
We keep the assumptions of section \ref{secprel}.
Not to repeat the same argument twice,
we treat the equivariant case from the very beginning;
the reader willing to restrict to the non-equivariant case should
just let $G$ be the trivial group (and skip \ref{IC} as containing no
non-empty statement).

{\it We change the notations.}
From now on  $Coh$, $QuasiCoh$ will denote  the category of $G$-equivariant
coherent (respectively, quasicoherent)
{\it equivariant} sheaves on $X$. Also $\Xtop$ will denote a subset
in the topological space of the scheme $X$, consisting 
of generic points of $G$-invariant subschemes;
we will endow $\Xtop$ with the induced topology.
Thus $\Xtop$ maps to the topological space of $S$, and for $s\in S$
the fiber over $s$ is the set of points of $X_s$ which are invariant under the 
component of identity in $G_s$.

 We will say that
 $x, y\in \Xtop$ are equivalent (and write $x\sim y$) if $x\in G(y)$
(i.e. if $x\in \Ztop \iff y\in \Ztop$ for a $G$-invariant subscheme 
$Z\subset X$).  The set of equivalence classes
$\Xtop/\sim$ is identified with the set of points of the stack
$X/G$. 

According to  Proposition \ref{dual_exist},
$X$ has an equivariant
 dualizing complex; we fix one, denote it by $\DC$.
 This choice defines
the {\it codimension function} $d$ on (all) points of $X$,
 which is determined by
the condition that $i_x^!(\DC)$ is concentrated in 
homological degree $d$  (see \cite{H}, \S V.7).
We set $\dim(x)=-d(x)$; if, say, $X$ is of finite type over a field, we
can (and will) assume that $\dim(x)$ is the (Krull) dimension of
the closure of $x$. Notice that $\dim(x)=\dim(y)$ for $x,y\in \Xtop$,
$x\sim y$.

Let $\tauleqn: D^?(Coh)\to D^{\leq n}(Coh)$, $\taugeqn:D^?(Coh)\to
D^{\geq n}(Coh)$ be the  truncation
 functors with respect to the usual $t$-structure on $D^?(Coh)$. (Here 
$?=+,-$ or $b$.)

Let  $p$ (perversity) be an
 integer-valued function on $\Xtop$, constant on equivalence classes.

 We define
the dual perversity by $\pbar(x )=-\dim(x  )-p(x )$.

\begin{Def} We define  $\Dpm \subset D^-(Coh)$,  $\Dpp \subset D^+_c(QuasiCoh)$
 by:

 $\F \in \Dpp$ if 
for any $x \in \Xtop$  we have $i_x ^!(\F )\in D^{\geq p(x )}(\O _x -mod)$.

$\F \in  \Dpm$ if 
for any $x \in \Xtop$  we have  $i_x ^*(\F )\in D^{\leq p(x )}(\O _x -mod)$.

\end{Def}

\begin{Lem}\label{predvar} a) $\D(\Dpm)=D^{\pbar,\geq 0}$.

b) Let  $i_Z:Z\imbed X$ be  a locally closed ($G$-invariant)
subscheme. Then $p$ defines
the induced perversity $p_Z=p\circ i_Z:Z^{top}\to \Zet$.
We have: 

$i_Z^*$ sends $\Dpm$ to $D^{p_Z,\leq 0}$; $i_Z^!$ sends 
$\Dpp$ to  $D^{p_Z,\geq 0}$; 

c) In the situation of (b) assume that $Z$ is closed. Then 
$i_{Z*}$ sends $D^{p_Z,\leq 0}$ to $\Dpm$, and  $D^{p_Z,\geq 0}$ to
$\Dpp$. 
\end{Lem}

\proof a) One knows from \cite{H}, \S V.6 that  for any $\F$ in the bounded
derived category of coherent sheaves we have $i_x^!(\D(\F))
= Hom_{\O_x}(i_x^*(\F), I_{\O_x})[-\dim(x)]$, where $I_{\O_x}$ is the
injective hull of the residue field of $\O_x$. Since $Hom_{\O_x}
(\underline{\ \ }, I_{\O_x})$ is exact and kills no finitely generated 
$\O_x$ module, (a) follows.

b) follows from Lemma \ref{equiv}; in view of this Lemma
if $\F\in D^{p,\leq 0}$, then for any $x \in \Ztop\subset \Xtop$
there exists a subscheme $Z'\subset Z$ with generic point $x $, such 
that $i_{Z'}^*(\F)=i_{Z'}^*(i_Z^*(\F))\in D^{\leq p(x )}(Coh_{Z'})$,
which implies $i_Z^*(\F)\in D^{p_Z,\leq 0}$; and similarly for $i_Z^!(\F)$.

c) is obvious. \epf

\begin{Prop}\label{Hom0} For $\F\in \Dpm$, $\G\in D^{p,>0}$ we have
  $Hom(\F , \G)=0$. 
\end{Prop}

\proof
We proceed by Noetherian induction in $X$; thus we can assume
that the statement with $(X,p)$ replaced by $(Z,p_Z)$ for a 
closed ($G$-invariant) subscheme $Z\subsetneq X$ is known.
(Otherwise replace $X$ by a minimal closed ($G$-invariant) subscheme
for which it is false).

Fix $\F \in \Dpm$, $\G\in  D^{p,>0}$. Let $x $ be a generic point of  $X$.
Using  Lemma \ref{equiv} we find an   open ($G$-invariant)
subscheme $j:U\imbed X$ containing $x $,
 such that $j^*(\F )\in D^{\leq
p(x )}(Coh(U))$ and $j^*(\G)=j^!(\G)\in D^{>p(x )}(Coh(U))$.
 Thus, of course, $Hom(j^*(\F ),
j^*(\G))=0$. 

Let $i$ denote the closed imbedding of $\Xtop-U^{top}$ into $\Xtop$. 
Consider the distinguished triangle in $D^b(QuasiCoh)$:
$i_*i^! \G \to \G \to j_*j^*(\G)\to i_*i^! \G[1] $. 
By Lemma \ref{lim}(a) we see that 
$$Hom (\F ,i_*i^! \G)=
\varinjlim_Z Hom(\F ,i_{Z*}i_Z^!(\G))=\varinjlim_Z
Hom(i^*_Z(\F ),i_Z^!(\G))=0,$$ because $i^*_Z(\F ) \in D^{p_Z, \leq 0}$,
$i_Z^!(\G) \in  D^{p_Z, > 0}$ by Lemma \ref{predvar}b), so  
$Hom(i^*_Z(\F ),i_Z^!(\G))=0$ by the induction hypotheses.
This implies the desired equality $Hom(\F ,\G)=0$, since
$Hom(\F ,j_*j^*(\G))=Hom(j^*(\F ),j^*(\G))=0$. \epf

\medskip

\begin{Def}\label{mon}
 A perversity function $p$ is 

 {\it monotone} if
$x '\in \overline x \; \Rightarrow p(x ')\geq p(x )$; 

 {\it strictly monotone} if
$x '\in \overline x \; \Rightarrow p(x ')> p(x )$; 

 (strictly) {\it comonotone} if the dual perversity
$\pbar(x )=-\dim(x )-p(x )$ is (strictly) comonotone.
\end{Def}

\begin{Thm}\label{t_str} Suppose that a perversity $p$ is monotone
and comonotone. Then 
$(\Dpm\cap D^b , \Dpp\cap D^b)$ define a $t$-structure on $D^b(Coh)$.
\end{Thm}

\proof
In view of Proposition \ref{Hom0} we have only to show, that for
any $\F\in D^b(Coh)$ there exists a distinguished triangle
$\F'\to \F\to \F''$ with $\F'\in D^{p,\leq 0}$, $\F''\in D^{p,>0}$.
We again
 proceed by Noetherian induction; thus we can assume 
that for a closed ($G$-invariant) subscheme $Z\subsetneq X$,
and $\F\in D^b(Coh_Z)$ there exists a triangle $\F'\to \F\to \F''$
with  $\F'\in D^{p_Z,\leq 0}(Coh_Z)$, $\F''\in D^{p_Z,>0}(Coh_Z)$.

It will be convenient to use the following notation
(see \cite{BBD}, 1.3.9).  If $D'$, $D''$ are sets of (isomorphism classes)
of objects of a triangulated category $D$, then $D'*D''$ is the set
of  (isomorphism classes) of objects of $D$, defined by the condition:
$B\in D'*D''$ iff there exists a distinguished triangle $A\to B\to C$
with $A\in D'$, $C\in D''$. 
The octahedron axiom implies (see \cite{BBD}, Lemma 1.3.10)
that the $*$ operation is associative, i.e. $(D'*D'')*D'''
=D'*(D''*D''')$. Thus the meaning of the
notation $D_1*\dots *D_n$ is unambiguous.

We will make the following abuse of notations: for a category
$\A$  we will write $\A$ instead of 
``the set of isomorphism classes of $Ob(\A)$''; and for a functor
$F:\A\to \B$ we will write $F(\A)$ instead ``image of the map
from the set of isomorphism classes of  $Ob(\A)$ to that of $Ob(\B)$
induced by $F$''. 
Then the statement we want to prove says that 
\begin{equation}\label{xotim}
D=D^{p,\leq 0}*D^{p,>0}.
\end{equation}

We claim that it is enough to show that 
\begin{equation}\label{reduce}
D=\cupl_Z D^{p,\leq 0}* i_{Z*}(D^b(Coh_Z)) * D^{p,>0},
\end{equation}
where $Z$ runs over all  ($G$-invariant) closed subschemes $Z\subsetneq X$.
Indeed, by the induction assumption we know that $ D^b(Coh_Z)
= D^{p_Z,\leq 0}(Coh_Z) * D^{p_Z,>0}(Coh_Z)$.
Thus \eqref{reduce} implies that
 $$D=\cupl_Z  D^{p,\leq 0}*i_{Z*}(D^{p_Z,\leq 0}(Coh_Z))
 *i_{Z*}( D^{p_Z,>0}(Coh_Z))
*  D^{p,>0}.$$
Rewriting the latter expression as
$$ \cupl_Z  \left( D^{p,\leq 0}*i_{Z*}(D^{p_Z,\leq 0}(Coh_Z)) \right) *
\left( i_{Z*}(D^{p_Z,>0}(Coh_Z)) *  D^{p,>0} \right),$$ 
and noting that by Lemma \ref{predvar}(c) we have
 $i_{Z*}(D^{p_Z,\leq 0}(Coh_Z))\subset D^{p,\leq 0}$, hence
$$ D^{p,\leq 0}*i_{Z*}(D^{p_Z,\leq 0}(Coh_Z))\subset D^{p,\leq 0}
* D^{p,\leq 0} = D^{p,\leq 0},$$
and simlarly for $D^{p,>0}$, we get  \eqref{xotim}.

Let us prove \eqref{reduce}. Fix $\F\in D$, and a generic point
$x $ of $X$. Let $\jbar:\overline{G(x )}\imbed X$ be the imbedding of the 
closure of $G(x )$ (in particular, if $X$ is irreducible, then $\jbar=id$).
We set $\F^-=\tau_{\leq p(x )}^{stand}(\jbar_\bullet \jbar ^!(\F))$.
Thus $\F^-\in D^{p,\leq 0}$, because it is supported on $\overline{G(x )}$,
and $p$ is monotone. Also
we have a canonical morphism $\F^-\to \F$. Let $\F_1$ be its cone;
then $i_x ^*(\F_1)\in D^{>p(x )}(\O_x -mod)$.

The dual procedure (in the sense of Grothendieck-Serre duality)
gives $\F^+ \in D^{p,>0}$, and a morphism $f:\F_1\to \F^+$, such that
$i_x ^*(f)$ is an isomorphism. More presicely, we set
$$\F^+= \D( \tau^{stand}_{< \pbar(x ) } \jbar_\bullet\jbar^!(\D(\F_1))).$$
Since $p$ is comonotone, we see by Lemma \ref{predvar}(a) that
$\F^+\in D^{p,>0}$. Since the local duality for the Artinian ring 
$\O_x $ is an exact functor (\cite{H}, \S V.6), we see that
$i_x ^*(\D(\F_1))\in D^{\leq \pbar(x )}(\O_x -mod)$, thus $i_x ^*(\D\F^+)
\iso i_x ^*(\D \F_1)$, and hence also $i_x ^*(\F_1)\iso i_x ^*(\F^+)$.

Thus, if we set $\F^0=cone(\F_1\to \F^+)[-1]$, then  $i_x ^*(\F^0)=0$.
Hence by Lemma \ref{lim}(b) we have
$\F^0\cong i_{Z*}(\F_Z)$  for some  closed ($G$-invariant)
subscheme $Z\subsetneq X$, and an object $\F_Z \in D^b(Coh_Z)$. 
So we  get $$\F\in \{\F^-\}*\{\F^0\}
*\{\F^+\}\subset D^{p,\leq 0}*i_{Z*}(D^b(Coh_Z)) * D^{p,>0},$$
 which proves \eqref{reduce}. \epf

\begin{Rem} Construction of an object $\F^+\in D^{p,>0}$
 with given generic fiber (and with a morphism from a given object)
is the only place in this paper,
 where  the (equivariant) duality formalism is used.
 
\end{Rem}

\begin{Cor}\label{cohext} \footnote{This statement and idea of proof are copied
from a message by Deligne to  the author. 
(Possible mistakes belong to the author).}
Let $j:U\imbed X$ be an open subscheme, $p:\Xtop\to \Zet$ be a monotone and
comonotone perversity, and $\F \in D^{p,\geq 0}(Coh(U))$.
Consider $j_*(\F)\in D^b(QuasiCoh(X))$, and let $n=\min\limits_{x\not\in \Utop}
p(x)$. Then $\tau^{stand}_{\leq n-2}(j_*(\F))$ has coherent cohomology.
\end{Cor}

\proof Let $\Ftil\in D^b(Coh)$ be any extension of $\F$
(see Corollary \ref{exte}); replacing
$\Ftil$ by $\tau_{\geq 0}^p
(\Ftil)$ we can achieve that $\Ftil\in D^{p,\geq 0}$.
If $\bfZ$ denotes $\Xtop-\Utop$, we can consider the exact triangle
$i_{\bfZ*}i_\bfZ^!(\Ftil) \to \Ftil \to j_*(\F)$. 
Since $ \Ftil$ has coherent cohomology, it is enough to show that
$\tau^{stand}_{\leq n-2}(i_{\bfZ*}i_\bfZ^!(\Ftil)[1])\in D^b(Coh)$ 
as well. However, the assumption $\Ftil\in D^{p,\geq 0}$ implies
that $i_Z^!(\Ftil)\in D^{p,\geq 0}(Z)\subset D^{\geq n}(Coh_Z)$
for $\Ztop=\bfZ$, hence $i_\bfZ^!\in D^{\geq n}(QuasiCoh)$, and
 $\tau^{stand}_{\leq n-2}(i_{\bfZ*}i_\bfZ^!(\Ftil)[1])=0$. \epf

\begin{Rem}\label{Grfin}
It was pointed out to us by Deligne that Corrolary \ref{cohext}
 is equivalent to
the Grothendieck Finiteness Theorem, \cite{SGA2},  VIII.2.1.
\end{Rem}

\subsection{Coherent IC-sheaves}\label{IC} 
We will assume that $p$ is a monotone and comonotone perversity function.
We will denote the core of the $t$-structure on $D^b (Coh(X))$
constructed in the previous section by $\P=\P_X=\P_{X,p}$.

\begin{Lem}\label{orto}
 Let $\bfZ=Z_0^{top}\subset \Xtop $ for a closed ($G$-invariant)
subscheme $Z_0\subset X$.
 Let $\F\in \P_X$.

a) The following  conditions are equivalent: 

\ i) $i_x ^*(\F)\in D^{<p(x )}(\O_x -mod)$ for all $x \in \bfZ$.

\ ii) $i_Z^*(\F)\in D^{p_Z,<0}$ for any closed ($G$-invariant)
subscheme $Z\subset X$,
 $\Ztop\subset \bfZ$.

\ iii) $Hom(\F, \G)=0$ for all $\G\in \P$, such that supp$(\G)\cap \Xtop
\subset \bfZ$;

b) The following  conditions are equivalent:

\ i) $i_x ^!(\F)\in D^{>p(x )}(\O_x -mod)$ for all $x \in \bfZ$.

\ ii) $i_Z^!(\F)\in D^{p_Z,>0}$ for any closed ($G$-invariant)
subscheme $Z\subset X$,
 $\Ztop\subset \bfZ$.

\ iii) $Hom(\G, \F)=0$ for all $\G\in \P$, such that supp$(\G)\cap \Xtop
\subset \bfZ$.
\end{Lem}
\proof  (a,i) $\iff$ (a,ii),
 follows from Lemma \ref{equiv} (a).
For a closed subscheme $Z$ and an object
$\G\in D^b(Coh_Z)$ we have $Hom(\F, i_{Z*}(\G))=
Hom(i_Z^*(\F),\G)$.
If $\F\in \P$, then $i_Z^*(\F)\in D^{p_Z,\leq 0}$; however,  for an
object of any triangulated category with a $t$-structure, and an object $A\in
D^{\leq 0}$ we have $A\in D^{<0}\iff Hom(A,B)=0$ for all $B$ in the core of the
 $t$-structure. This shows (a,ii) $\iff$ (a,iii).
Thus (a) is proved, and the proof of (b) is similar. \epf

\medskip

It is  convenient to reformulate the conditions of Lemma \ref{orto} as
follows. Let the auxilary perversity functions $p^-=p^-_{(\bfZ)}$, $p^+
=p^+_{(\bfZ)}$ be given by $p^-(x )
=p(x )=p^+(x )$ if $x \not\in \bfZ$, and $p^-(x )=p(x )-1$, $p^+(x )=p(x )+1$ if $x 
\in \bfZ$. Then conditions (a) of Lemma \ref{orto} just say that
$\F\in D^{p^-,\leq 0}$, and conditions (b) say that $\F\in D^{p^+,>0}$.

\begin{Thm}\label{GM}
 Let $j:U\imbed X$ be a ($G$-invariant) locally closed subscheme,
set $p^-=p^-_{(\Ubar-U)}$, $p^+=p^+_{(\Ubar-U)}$
and define a full subcategory $\P_{!*}(U)\subset \P_\Ubar$ by
  $\P_{!*}(U)= D^{p^-,\leq 0}(Coh_\Ubar) \cap  D^{p^+,\geq 0}(Coh_\Ubar)$.

 Suppose that $p(x ^-)>p(x )$,
  $\pbar (x ')>\pbar (x )$ for any $x \in \Utop$, 
$x '\in \xbar\cap \Ubar^{top}$,
 $x '\not \in \Utop$. 
Then $j^*$ induces an equivalence between $\P_{!*}(\U)$ and $\P_U$.

The inverse equivalence is denoted by $j_{!*}:\P_U\to \P_{!*}(U)
\subset \P_\Ubar$, and is called the functor of minimal (or Goresky-MacPherson,
or IC) extension.
\end{Thm}
\proof The conditions of the Theorem say that both $p^-$ and $p^+$
induce monotone and comonotone perversity functions on $\Ubar^{top}$;
hence they  define  $t$-structures on $D^b(Coh_\Ubar)$.
 Let $\tau^-$,
$\tau^+$ be the corresponding truncation functors.

We first introduce an auxilary functor $J_{!*}$ on $D^b(Coh_\Ubar)$
by $J_{!*}=\tau^-_{\leq 0} \circ \tau^+_{\geq 0}$.

\begin{Lem}\label{JGM} a) $J_{!*}$ takes values in $\P_{!*}(U)$.

b) If $f:\F\to \G$ is a morphism in $D^b(Coh_\Ubar)$, 
such that $j^*(f)$ is an isomorphism, then $J_{!*}(f)$ is an isomorphism.
\end{Lem}

\proof It is obvious that $J_{!*}(\F)\in D^{p^-,\leq 0}$; also, if $\F_1$
denotes $ \tau^+_{\geq 0}(\F_1)$, then we have an the exact triangle
$\tau^-_{>0}(\F_1) [-1] \to J_{!*}(\F) \to \F_1$. Here certainly
$\F_1\in D^{p^+,\geq 0}$; and also $\tau^-_{>0}(\F_1) [-1]\in D^{p^-, \geq 2}
\subset D^{p^+, \geq 0}$. This proves (a).

If $f$ is as in (b), then also $j^*(J_{!*}(f))$ is an isomorphism.
But then $J_{!*}(f)$ is a morphism in $\P_{!*}(U)$, such that
$j^*(f)$ is an isomorphism; thus its kernel and cokernel are objects
of $\P_\Ubar$ supported on $\Ubar-U$. However, Lemma \ref{orto}
says that $J_{!*}(\F),$
$J_{!*}(\G)$ have no subobjects or quotients supported on $\Ubar -U$. \epf

\medskip

Now, using Corollary \ref{exte}, we see from Lemma \ref{JGM} that there
exists a canonically defined functor 
$\jtil_{!*}:D^b(Coh_U)\to \P_{!*}(U)$, equipped
with an isomorphism $\jtil_{!*}\circ j^*=J_{!*}$. We set $j_{!*}=
\jtil_{!*}|_{\P_U}$. Then it is clear that $j^*\circ j_{!*} =id_{\P_U}$
canonically. Also $j_{!*}\circ j^*|_{\P_{!*}(U)} =id$ canonically, because
$J_{!*}|_{\P_{!*}(U)} =id$. Thus $j^*$ and $j_{!*}$ are  inverse
equivalences between $\P_{!*}(U)$ and $\P_U$. \epf

From now on assume that  $S=Spec(k)$, where $k$ is  a field.
For a $G$-orbit $O\subset X$ we set $p(O)=p(x )$, where $x $ is a generic
point of $O$ (this number does not depend on the choice of $x $
 because $x\sim x'$ if $x,x'$ are generic points of $O$).

\begin{Cor} \label{irre} 
For $\F \in D^b(Coh^G)$ the following statements are equivalent:

i) $\F $ is an irreducible object of $\P^G$.

ii) There exists a $G$-orbit $j:O\imbed X$,
such that 
$p(O)<p(x )$, $\pbar(O)<\pbar(x )$ for any non-generic point 
$x \in \Obar^{top}$,
and an irreducible $G$-equivariant vector bundle $L$ on $O$, such that
$\F=j_{!*}(L[p(O)])$.
\end{Cor}
\proof  (ii) $\Rightarrow$ (i) is obvious from Lemma \ref{orto}; let us
prove
(i) $\Rightarrow$ (ii). Let $x$ be a generic point 
of $supp(\F)$, and  $Z\subset supp(\F)$, $Z\not \owns x$ be a 
closed ($G$-invariant) subscheme.
If $\F$ is irreducible, then  for $\G\in \P_Z$ 
we have
 $Hom(\F, i_{Z*}(\G))=0$,
$Hom( i_{Z*}(\G), \F)=0$. 
Thus Lemma \ref{orto} says that  $i_Z^*(\F)\in D^{p_Z,<0}$, 
$i_Z^!(\F) \in  D^{p_Z,> 0}$; and $i_x ^*(\F)\in D^{<p(x )}(\O_x -mod)$,
$i_x ^!(\F)\in D^{>p(x )}(\O_x -mod)$. In particular, it follows
that $Z$ can not contain generic points of $supp(\F)$, hence $supp(\F)=
\overline{G(x)}$. We also see that if $x'$ is a non-generic point of 
$supp(\F)$, then  $p(x )<p(x ')$;
 indeed, otherwise
the coherent sheaf ${\cal H}^{p(x )}(\F)$ has a nonzero fiber at $x $,
but has  zero fiber at $x '$, which contradicts the Nakayama Lemma.
Applying the Grothendieck-Serre duality we get also $\pbar(x )<\pbar(x ')$.
In particular, for any non-generic point
 point $x'\in (supp (\F))^{top}$ we have
$\dim(x')< \dim(x)-1= \dim  (supp (\F))-1$. Then Rosenlicht's
Theorem (see e.g. \cite{Rosenlicht}) implies that
 $x$ is a generic point of an orbit $O$.
 Thus $\F\in \P_{!*}(O)$, 
so (ii) follows from Theorem \ref{GM}.
\epf

\begin{Cor}\label{Art} Suppose that  
$G$ acts on $X$ with a finite number of orbits, and
 $p$ is strictly monotone and comonotone.
Then the  category $\P^G$ is Artinian.
\end{Cor}

\proof Conditions of the Corollary imply that $j_{!*}^O$ is defined
for any $G$-orbit $j^O:O\imbed X$. By induction in the number of orbits
one can deduce 
(using Corollary \ref{exte}) that the irreducible
objects $j^O_{!*}(L[p(O)])\in \P$  generate the triangulated category
 $D^b(Coh)$.
This implies that $\P$ is Artinian.
\epf

\begin{Ex}\label{N}
 Let $G$ be a simple group over a field of characteristic 0
(or of large finite characteristic), and let $\N\subset G$ be the subvariety of
unipotent elements. Then $G$ acts on $\N$ by conjugation, and this action has
a finite number of orbits. Moreover, dimension of an orbit is known to be even.
Thus the set $\N^{top}$ consists of generic points of 
$G$-orbits, and  we can define the ``middle perversity'' by
 $p(x _O)=-\frac{\dim(O)}{2}$ 
for an orbit $O\subset \N$
 (where $x _O$ is the generic point of $O$). Then $p$ is
obviously strictly monotone and comonotone, hence by Proposition \ref{Art}
the kernel of the corresponding $t$-structure is Artinian.
See \cite{me} for more information on this example.
\end{Ex}

\begin{Rem} It will be shown in \cite{AB} that the irreducible objects of the 
$t$-structure described in Example \ref{N} 
are closely related to cohomology
of (tilting) modules over a quantum group at a root of unity. 
(This relation was independently conjectured by Ostrik).
\end{Rem}

\end{document}